\documentclass[11pt]{article}
\usepackage{color, amsmath,amssymb, amsfonts, amstext,amsthm, latexsym}

\usepackage{amssymb, epsfig, amssymb, latexsym}
\begin{document}

\newcommand{\ddt}{\frac{d}{dt}}

\newcommand{\s}{\sigma}
\renewcommand{\k}{\kappa}
\newcommand{\p}{\partial}
\newcommand{\D}{\Delta}
\newcommand{\om}{\omega}
\newcommand{\Om}{\Omega}
\renewcommand{\phi}{\varphi}
\newcommand{\e}{\epsilon}
\renewcommand{\a}{\alpha}
\renewcommand{\b}{\beta}
\newcommand{\N}{{\mathbb N}}
\newcommand{\R}{{\mathbb R}}
   \newcommand{\eps}{\varepsilon}
   \newcommand{\Wsob}{\smash{{\stackrel{\circ}{W}}}_2^1(D)}
   \newcommand{\EX}{{\mathbb E}}
   \newcommand{\PX}{{\mathbb P}}

\newcommand{\cF}{{\cal F}}
\newcommand{\cG}{{\cal G}}
\newcommand{\cD}{{\cal D}}
\newcommand{\cO}{{\cal O}}

\newtheorem{theorem}{Theorem}
\newtheorem{lemma}{Lemma}
\newtheorem{remark}{Remark}

\title{Stochastic Modeling of Unresolved Scales\\ in Complex Systems
\footnote{This research was partly supported by the NSF Grant
   0620539,   the Cheung Kong Scholars Program
   and the K. C. Wong Education Foundation.} }

\author{Jinqiao Duan \\
Department of Applied Mathematics\\ Illinois Institute of Technology \\
  Chicago, IL 60616, USA  \\\emph{E-mail: duan@iit.edu} \\
  and \\
School of Mathematics and Statistics\\ Huazhong University of Science and Technology \\
  Wuhan 430074, China  }

\date{October 20, 2008 (Revised version)    }

\maketitle

\begin{abstract}
Model uncertainties or simulation uncertainties  occur in
mathematical modeling of multiscale complex systems, since some
mechanisms or scales are not represented (i.e., ``unresolved") due
to lack in our understanding of these mechanisms or   limitations
in computational power. The impact of these unresolved scales on
the resolved scales needs to be parameterized or taken into
account. A stochastic   scheme is devised to take the effects of
unresolved scales into account, in the context of solving
nonlinear partial differential equations. An example is presented
to demonstrate this strategy.

\bigskip

%{\bf Short Title:}  Stochastic Modleing of Unresolved Scales\\


\medskip

%{\bf PACS Numbers:}  02.50.Fz, 05.10.Gg, 02.60.Lj, 02.60.Nm,
%02.60.Cb

\end{abstract}

\maketitle

\newpage

\section{Introduction}

 Mathematical
  models for scientific and engineering
  systems often involve with some uncertainties.
  We may roughly classify such uncertainties into two kinds.
  The first kind of uncertainties may be called \emph{model uncertainty}.
  They involve with  physical processes that are less known, not
  yet well understood,   not well-observed or measured, and thus
  difficult to be represented in the mathematical models.

The second kind of uncertainties may be called \emph{simulation
uncertainty}. This arises in numerical simulations of multiscale
systems that display a wide range of spatial and temporal scales,
with no clear scale separation. Due to the limitations of computer
power, at present and for the conceivable future,  not all scales
of variability can be explicitly simulated or resolved. Although
these unresolved scales  may be very small or very fast, their
long time impact on the resolved simulation may be delicate (i.e.,
may be negligible or may have significant effects, or in other
words, uncertain). Thus, to take the effects of unresolved scales
on the resolved scales into account,   representations or
parameterizations of these effects   are desirable.

These uncertainties are sometimes also called \emph{unresolved
scales}, as they are not represented or not  resolved  in modeling
or simulation.
  Model uncertainties have been considered in, for example,
  \cite{Gar, Horst, Hasselmann, Arn00, WaymireDuan, Lin, Palmer2, Pas, Sura}
  and references therein.  Works relevant for parameterizing unresolved scales
  include \cite{Kantz, Stuart, Majda,  Sardes, Berloff, DuanBalu, Wilks,
Williams, Sagaut, Berselli}, among others.

In this paper we consider an issue of  approximating model
uncertainty or simulation uncertainty (unresolved scales)   by
stochastic processes, and then devise a stochastic scheme for such
approximations. We first recall some basic facts about fractional
Brownian motion (fBM) in \S \ref{fbm}. Then we discuss  model
uncertainty and simulation uncertainty in \S \ref{modelerror} and
\S \ref{SLES}, respectively. Finally, we present an example in \S
\ref{example} demonstrating our result. This example involves
approximating subgrid scales via correlated noises, in the context
of large eddy simulations of a partial differential equation.

\section{Fractional Brownian motion and colored noise }
\label{fbm}

We   discuss a model of colored noise in terms of fractional
Brownian motion (fBM), including a special  case which is white
noise in terms of usual Brownian motion.
 The fractional Brownian motion $B^H(t)$, indexed by a so
called Hurst parameter $H \in (0, 1)$,  is a generalization of the
more well-known process of the usual Brownian motion $B(t)$. It is
a centered Gaussian process with stationary increments. However,
the increments of the fractional Brownian motion are not
independent, except in the usual Brownian motion case
($H=\frac12$). For more details, see \cite{Nualart, Memin,
Duncan, Maslowski, Tindel}.\\
\\
 Definition of fractional Brownian motion: For $H \in (0,
1)$, a Gaussian process $B^H(t)$, or $fBM(t)$, is a fractional
Brownian motion if it starts at zero  $B^H(0)=0,\; a.s.$, has mean
zero  $\EX[B^H(t)] = 0 $, and has covariance $\EX[B^H(t)B^H(s)] =
\frac12(|t|^{2H} + |s|^{2H} - |t-s|^{2H})$ for all t and s. The
standard Brownian motion is a fractional Brownian motion with
Hurst parameter $H=\frac12$. \\
\\
 Some properties of fractional Brownian motion:
A fractional Brownian motion
 $B^{H}(t)$ has the following properties:\\
(i) It has stationary increments;\\
(ii) When $H=1/2$, it has independent increments; \\
(iii) When $H \neq 1/2$,   it is neither Markovian, nor a
semimartingale. \\
%(iv) The covariance between future and past increments is positive
%if $H>1/2$ and negative if $H<1/2$.\\

%\textbf{Ito Isometry with respect to fractional Brownian motion }\\
%For $ H >1/2$,    $f \in L^2$ and $f(t)$ being deterministic, the
%following isometry property holds:
%\begin{eqnarray*}
%\EX(\int^T_0 f(u)dB^{H}_t\int^T_0
%g(u)dB^{H}_t)=H(2H-1)\int^T_0\int^T_0 f(u)g(v)|u-v|^{2H-2}dudv
%\end{eqnarray*}
%and
%\begin{eqnarray*} \EX(\int^T_0
%f(t)dB^{H}_t)^2=H(2H-1)\int^T_0\int^T_0 f(u)f(v)|u-v|^{2H-2}dudv.
%\end{eqnarray*}

%When the hurst parameter $H
%> 1/2$, it has positive correlation. This property is suitable to
%simulate the subgrid scale term $R(x,t)$ in our system, because  it is
%correlation in time. Please see the correlation plot below:

%\noindent \textbf{How to simulate a fractional Brownian motion}\\
%The increments of fractional Brownian motion are correlated.

We use the Weierstrass-Mandelbrot function to approximate the
fractional Brownian motion. The basic idea is to simulate
fractional Brownian motion by randomizing a representation  due to
Weierstrass. Given the Hurst parameter $H$ with  $0<H<1$, we
define the function $w(t)$ to approximate the fractional Brownian
motion:
\begin{eqnarray*}
w(t_i) = \sum^{\infty}_{j=-\infty} C_j r^{jH}\sin(2\pi r^{-j}t_i +
d_j)
\end{eqnarray*}
where $r=0.9$ is a constant,  $C_j$'s are normally distributed
random variables with mean $0$ and standard deviation $1$, and the
$d_j$'s are uniformly distributed random variables  in the
interval $0 \leq d_j < 2 \pi$. The underlying theoretical
foundation for this approximation can be found in \cite{Pipiras,
Mehrabi}. Figures 1 and 2   show a sample path of the usual
Brownian motion (i.e., $H=\frac12$), and fractional Brownian
motion with Hurst parameter $H=\frac34$, respectively.

\begin{figure}[htbp]
\begin{center}
\includegraphics[height=3in,width=4in]{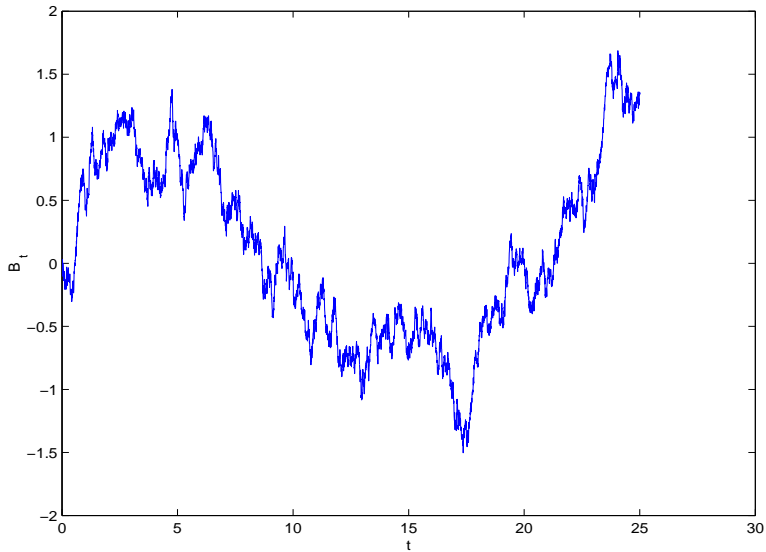}
\caption{\small \sl A sample path of   Brownian motion $B(t)$}
 \end{center}
\end{figure}

\begin{figure}[htbp]
\begin{center}
\includegraphics[height=3in,width=4in]{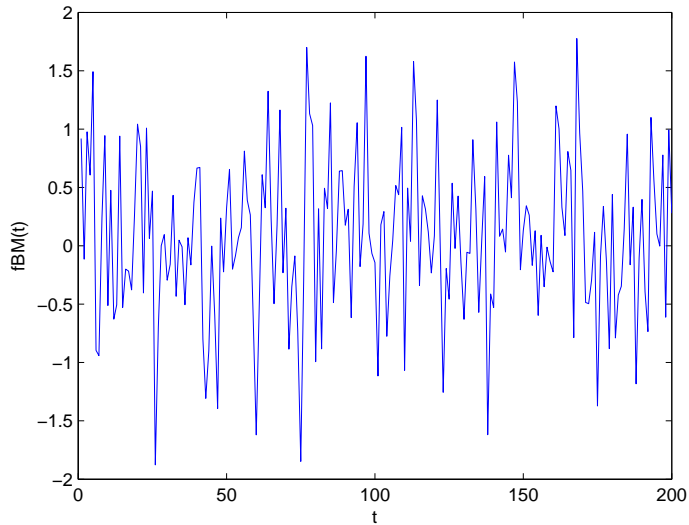}
\caption{\small \sl A sample path of fractional Brownian motion
$B^H(t)$, with $H=0.75$}
\end{center}
\end{figure}

%\begin{figure}[htbp]
%\begin{center}
%\includegraphics[height=3in,width=4in]{fbm-2.eps}
%\caption{\small \sl Another sample path of Fractional Brownian
%motion, with $H=0.75$.}
%\end{center}
%\end{figure}

%%%%%%%%%%%%%%%%%%%%%%%%%%%%%%%%%%%%%%%%%%%%%%%%%%%%%%%%%%%%%%%%%%%%%%%%%%%%%
\section{Model uncertainty } \label{modelerror}

We consider a   spatially extended system modeled by a partial
differential equation (PDE):
\begin{eqnarray} \label{det}
 u_t =  Au+ N(u),
\end{eqnarray}
where $A$ is a linear (unbounded) differential operator, and $N$
is a nonlinear function of $u(x,t)$ with $x\in D$ and $t>0$, and
satisfies a local Lipschitz condition. In fact, $N$ may also
depend on the gradient of $u$.

If this (deterministic) model is accurate, i.e., its prediction on
the field  $u$ matches with the observational data $\tilde{u}$ on
a certain period of time $[0, T]$, then there is no need for a
stochastic approach. However, when the prediction $u$ deviates
from the observational data $\tilde{u}$, we then need to modify
the model \eqref{det}. In this case, the observational data
$\tilde{u}$ may be thought to satisfy a modified model:
\begin{eqnarray} \label{stoch1}
 \tilde{u}_t =  A\tilde{u}+ N(\tilde{u}) + F(\tilde{u}),
\end{eqnarray}
where the model uncertainty $F(\tilde{u})$ is usually a
fluctuating (i.e., random) process, as the observational data
$\tilde{u}$ is so (i.e., has various samples or realizations).

The model discrepancy or model uncertainty $F(\tilde{u})$ may have
various causes, such as missing physical mechanisms (not
represented in the deterministic model \eqref{det}). Sometimes,
the model uncertainty $F(\tilde{u})$ is smaller in magnitude than
other terms in the model \eqref{stoch1} and thus is often ignored
in the deterministic modeling. However, being small and being
fluctuating may not necessarily imply that its impact on the
overall system evolution to be small \cite{Arnold}. To take this
impact into account, we would like to model or approximate
$F(\tilde{u})$ by a stochastic process.

We first calculate the model uncertainty $F(\tilde{u})$ via
observational data $\tilde{u}$. By discretizing  \eqref{stoch1}
and using data samples for $\tilde{u}$, we obtain (discretized)
samples for $F$.

The time correlation may then be calculated using the samples of
$F$. If the time correlation scale is significantly shorter than
the time scale for the field $u$, we may ignore the time
correlation and thus approximate $F$ by the following stochastic
process containing a (uncorrelated) white noise, for example:
\begin{eqnarray} \label{stoch2}
 F  = f  + \sigma \; \tilde{u} \dot{B}_t,
\end{eqnarray}
where $f=\EX F$ is the mean of $F$ (computed from data), $B_t$ is
the usual Brownian motion (reviewed in \S \ref{example} below) and
the deterministic noise intensity $\sigma$ may depend on space.
Here $\sigma$ may be computed via stochastic calculus, especially
the Ito isometry, as follows.
\begin{eqnarray*}
 F  - \EX F &=& \sigma \; \tilde{u} \dot{B}_t,\\
\EX [\int_0^T (F  - \EX F)  dt]^2 &=&  \sigma^2 \;\EX [\int_0^T \tilde{u} dB_t]^2,\\
\EX [\int_0^T (F  - \EX F)  dt]^2 &=&  \sigma^2 \;\EX  \int_0^T \tilde{u}^2(x, t) dt,\\
\end{eqnarray*}
Thus we obtain
\begin{eqnarray} \label{sigma}
  \sigma =  \sqrt{\frac{\EX [\int_0^T (F  - \EX F)  dt]^2}{\EX  \int_0^T \tilde{u}^2(x, t) dt} } \; .
\end{eqnarray}
With this approximation, we obtain the following stochastic
partial differential equation (SPDE) as a modified model for the
original deterministic model \eqref{det}: $U \approx \tilde{u}$
\begin{eqnarray} \label{spde}
 U_t =  AU+ N(U) +  f(U)  + \sigma \; \tilde{u} \dot{B}_t.
\end{eqnarray}

In general, the model uncertainty $F$ may    be better
approximated by correlated noise via fractional Brownian motion.
Since the procedure is similar, we will demonstrate this in the
next section when we discuss simulation uncertainty.

%%%%%%%%%%%%%%%%%%%%%%%%%%%%%%%%%%%%%%%%%%%%%%%%%%%%%%%%%%%%%%%%%%%%%%%%%%%%%
\section{Simulation uncertainty } \label{SLES}

This   section deals with simulation uncertainty, i.e.,
stochastically parameterizing the effects of the unresolved scales
on the resolved scales. We consider this issue in the context of
large eddy simulations (LES) of a  nonlinear partial differential
equation with memory.

 In large eddy simulations of fluid or geophysical
 fluid flows \cite{Sagaut, Berselli},
   the unresolved scales appear   as the so-called subgrid
  scales (SGS). The SGS term appears to be highly fluctuating
  (``random"); see the Figure 1 in \cite{Menev}. Partially motivated by   this,
 stochastic parameterizations of subgrid scales
  have been investigated in fluid, geophysical and climate simulations,   based on
  physical or intuitive or empirical arguments. Another,
  perhaps more important, motivation for applying stochastic
  parameterizations of subgrid scales is to induce the desired
  backward energy flux (``stochastic backscatter")
  in fluid simulations \cite{Leith, Mason, Schumann}.

We present one stochastic parameterization scheme of the subgrid
scale term in the   large eddy simulation of a nonlinear partial
differential  equation with an extra memory term, which is in fact
a nonlinear integro-partial differential equation. The
approximation scheme  is based on stochastic calculus involving
with a fractional Brownian Motion, and the ``parameter' to be
calculated is a spatial function, which is derived using Ito
stochastic calculus.

\begin{eqnarray}
 u_t =  Au+ F(u),
\end{eqnarray}
where $A$ is a linear differential operator, and $F$ is a
nonlinear function of $u(x,t)$ with $x\in D$and $t>0$, and
satisfies a local Lipschitz condition. We investigate stochastic
parameterizations of unresolved scales  in the context of large
eddy simulations of the above system.

The idea of large eddy simulation is to split the flow into a
local, spatial mean (or average) and a fluctuation about the that
mean. The mean $\bar{u}$ is defined by filtering or mollification
(convolution with an approximate identity). The goal is to predict
the mean accurately. This is widely believed possible based on the
idea that since fluctuations have random character, their average
effects on the mean notion can successfully be medelled.

To filter the solution, we   pick a filter. Many different ones
are commonly used. To fix the ideas, in this paper, we use
Gaussian filter as in \cite{Berselli}, $ G_{\delta}(x)=
\frac{1}{\pi\delta^2} e^{-\frac{
    x^2}{\delta^2}}, $ where $\delta>0$ is the filter size and
the filter is such that: (i) $u* G_{\delta}$ is infinitely
differentiable in space and, (ii) $u*G_{\delta} \to u$ as $\delta
\to 0$ in $L^2(D)$. Here and hereafter $ u* G_{\delta}=\int_D
u(y,t)G_\delta(x-y)dy$ or the over bar $\bar{u}$  denotes
convolution.

\begin{remark}
The mean $\bar{u}$ is a weighted average of $u$ about the point
$x$. As $\delta \to 0$, the points near $x$ are weighted more and
more heavily, so $u*G_{\delta} \to u$ as $\delta \to 0$ in
$L^2(D)$.
\end{remark}

Using the fact that convolution commutes with differentiatian, we
get the space-filtered system:
\begin{eqnarray*}
 \bar{u}_t =  A\bar{u}+  \overline{F(u)},
\end{eqnarray*}
or
\begin{eqnarray} \label{LES2}
 \bar{u}_t =  A\bar{u}+ F(\bar{u}) + R(u,\bar{u}),
\end{eqnarray}
where  the subgrid scale term $ R(u,\bar{u}) :=
\overline{F(u)}-F(\bar{u})$. Since generally $ \overline{F(u)}
\neq F(\bar{u})$, the usual parameterization or closure problem of
the large eddy simulation has arisen. Due to inaccurate
(uncertain) initial conditions or boundary conditions,
$R(u,\bar{u}) $  is a correlated fluctuating process \cite{Menev,
DuanBalu}, depending on samples $\omega$ in a suitable sample
space $\Omega$. We thus would like to approximate the subgrid
scale term $ R(u,\bar{u})$ by a stochastic process with a
correlated (i.e., colored) noise component, for example:
\begin{eqnarray}
\label{noise2}
 R  = f(\bar{u}) + \sigma(x) \frac{d B^H_t}{dt},
\end{eqnarray}
  where $\frac{d B^H_t}{dt}$ is a  colored noise
  (generalized time derivative of a fractional Brownian motion;
  reviewed in \S \ref{example} below), and
\begin{eqnarray} \label{force2}
  f(\bar{u})=\EX R ,
\end{eqnarray}
   is the mean component
  of the subgrid scale term $R $. Moreover,
    the noise intensity $\sigma(x)$ is a  non-negative  deterministic function
to be determined from fluctuating SGS data $R$. The subgrid scale
term  $R $ may be inferred from observational data (see
\cite{Peters1, Peters2} for relevant information for subgrid
scales in Navier-Stokes equations), or from fine mesh simulations.

  %With this  stochastically parameterized SGS term $R$, the large eddy
  %solution $\bar{u}$ from (\ref{LES}) is a stochastic process as well.

  %Note that  $\bar{u}$, as a solution of a \emph{nonlinear} random
  %partial differential equation, is generally non-Gaussian. The
  %product $\sigma(\bar{u}, x,t) \dot{B}^H_t$ is generally also
  %non-Gaussian. Overall, $\sigma(\bar{u}, x,t) \dot{B}^H_t$
  %may approximate quite complicated random data $R$.

Note that $\sigma$ is to be calculated or estimated from the
fluctuating SGS data $R$, either from observation or  from fine
mesh simulations. So this is an inverse problem. As in usual
inverse problems \cite{Tar}, the stochastic parameterizations for
the SGS term $R$ is not unique. What we proposed above is merely
an example. This offers an opportunity for trying various
stochastic parameterization schemes, much as one uses various
smoother functions (e.g., polynomials or Fourier series) to
approximate less regular functions or data in deterministic
approximation theory.

To estimate the unknown parameter (function) $\sigma(x)$, we start
with    the following relation:
\begin{eqnarray}
 R  - \EX R  &=& \sigma(x) \frac{d B^H_t}{dt}.
\end{eqnarray}
 Taking  time integral over a computational interval $[0, T]$
  on  both  sides, we obtain
\begin{eqnarray*}
  \int_0^T[R -\EX R ]dt  =  \int_0^T \sigma(x)
 d{B}^H_t = \sigma(x) B^H_T.
\end{eqnarray*}
Therefore, taking mean-square on both sides,
 \begin{eqnarray*}
  \EX( \int_0^T [R -\EX R ]dt)^2 &=& \sigma^2(x)  T^{2H}.
\end{eqnarray*}
%Applying the Ito  isometry for stochastic integral with respect to
%fractional Brownian motion with $H > \frac12$, and using the fact
%that the integrand is deterministic, we obtain that
%\begin{eqnarray*}
 % \EX( \int_0^T [R(x,t)-\EX R(x,t)]dt)^2 &=&
% \sigma^2(x)\;H(2H-1)\int_0^T  \int_0^T |u-v|^{2H-2}dudv.
%\end{eqnarray*}
Thus an  estimator  for $\sigma(x)$ is
\begin{eqnarray}
 \sigma(x)= \frac1{T^H}\sqrt{\EX( \int_0^T [R -\EX R ]dt)^2} \;\;,   \label{ax2}
\end{eqnarray}
which can be computed numerically.

%With the calculated noise intensity term $\sigma(x)$ from
%(\ref{ax}), the stochastic LES model (\ref{new}) is established.

By the stochastic   parameterization  (\ref{noise2}) on the SGS
term $R$, with  $f$ determined from (\ref{force2}) and $\sigma $
from (\ref{ax2}), the LES model (\ref{LES2}) becomes a stochastic
partial differential equation (SPDE) for the large eddy solution
$U \approx \bar{u}  $:
\begin{eqnarray} \label{new2}
  U_t = A U +  N(U)
   + f(U) + \sigma(x)   \frac{d B^H_t}{dt},
\end{eqnarray}
with the appropriately filtered boundary condition   and filtered
initial condition.

%%%%%%%%%%%%%%%%%%%%%%%%%%%%%%%%%%%%%%%%%%%%%%%%%%%%%%%%
\section{An example}   \label{example}

We present a specific example of stochastic modeling of simulation
uncertainty of subgrid scales, in the context of large eddy
simulations. We consider the following nonlinear partial
differential equation with a memory term  (time-integral term)
\cite{DuDuan}:
\begin{eqnarray} \label{heat}
 u_t  &=& u_{xx} + u - u^3 + \int^t_0 {1 \over {1+|t-s|^{\beta}}} \;u(x,s)ds,
\end{eqnarray}
under appropriate initial condition $u(x, 0)=u_0(x)$ and boundary
conditions $u(-1,t)=a, \; u(1,t)=b$ with $a, b$ constants, on a
bounded domain $D : -1 \leq x \leq 1$. Here $\beta $ is a positive
constant. This model arises in mathematical modeling in ecology
\cite{Wu}, heat conduction in certain materials \cite{GMP,MK} and
materials science \cite{FS, MK}. The time-integral term here
represents a memory effect depending on the past history of the
system state, and this memory effect decays polynomially fast in
time.

The large eddy solution  $ \bar{u} $ is the true solution $u$
looked through a filter: i.e.,    through convolution with a
spatial filter $G_{\delta}(x)$, with spatial scale (or filter size
or cut-off size) $\delta>0$:
\begin{eqnarray*}
\bar{u} (x,t) : = u* G_{\delta}=\int_D u(y, t)G_\delta (x-y)dy.
\end{eqnarray*}
In this paper, we use a
 Gaussian filter as in \cite{Berselli}, $ G_{\delta}(x)=
 \frac{1}{\pi\delta^2} e^{-\frac{x^2}{\delta^2}}$.

On convolving (\ref{heat}) with $G_{\delta}$, the large eddy
solution $\bar{u}$ is   to satisfy
\begin{eqnarray*}
 \bar{u}_t &=& \bar{u}_{xx} +    \bar{u}
 - \overline{u^3} + \int^t_0 {1 \over {1+|t-s|^{\beta}}} \bar{u}(x,s)ds,
\end{eqnarray*}
 or
\begin{eqnarray}  \label{LES}
\bar{u}_t &=& \bar{u}_{xx} +  \bar{u}
 - {\bar{u}}^3 + \int^t_0 {1 \over {1+|t-s|^{\beta}}} \;\bar{u}(x,s)ds +R(x,t),
\end{eqnarray}
where  the remainder term, i.e., the subgrid scale (SGS) term
$R(x, t)$ is defined as
\begin{eqnarray}  \label{EF}
 R(x,t):=({\bar{u}})^3 - \overline{(u^3)} .
\end{eqnarray}

%\begin{remark}
%The integral term here is a grid scale term, since
%\begin{eqnarray*}
%[\int^t_0 k(t,s)u(x,s)ds] * G_{\delta} = \int^t_0 k(t,s)[u(x,s)*
%G_{\delta}]ds.
%\end{eqnarray*}
%\end{remark}

We can write $u= \bar{u} + u'$ with $\bar{u}$ the large eddy term
and $u'$ the fluctuating term. Note that $\bar{u}=u-u'$. So the
SGS term $R(x, t)$ involves  nonlinear interactions of
fluctuations $u'$ and the large eddy flows. Thus $R(x, t)$ may be
regarded as a function of $\bar{u}$ and $u'$: $R:=R(\bar{u}, u')$.

The leads to a possibility of approximating $R(x, t)$ by a
suitable stochastic process defined on a probability space $(\Om,
\mathcal{F}, \mathbb{P})$, with
 $\om \in \Om$, the sample space, $\sigma-$field $\mathcal{F}$
  and probability
measure $ \mathbb{P} $.  This means that we treat $R$ data as
random data as in \cite{Menev}, which take different realizations,
e.g., due to fluctuating observations or due to numerical
simulation with initial and boundary conditions with small
fluctuations. In fluid or geophysical fluid simulations, the SGS
term may be highly fluctuating and time-correlated \cite{Menev},
and this term may be  inferred from observational data
\cite{Peters1, Peters2}, or from fine mesh simulations.

%Before we devise how we parameterize the subgrid scale term $R(x,
%t)$ as a colored noise or time-correlated term, we first discuss
%fractional Brownian motion and colored noise in the next
%subsection.

%\begin{remark} $\delta$ is implicit in the process of modelling
%(statistical regression). When $\delta \to 0$ , $||\sigma(x,t)||
%\to 0\ in\ L^2(D)$. The extreme case is that when $\delta = 0$,
%$\sigma(x,t)\equiv 0$ because $R \equiv 0$.
%\end{remark}

This further suggests    for  parameterizing the subgrid scale
term $R(x, t)$ as a time-correlated or colored noisy term. The
increments of fractional Brownian motion are correlated in time
and hence its generalized time derivative $\dot{B}_t^H$ is used as
a model for colored noise. In the special case $H=\frac12$, we
have the white noise $\dot{B}_t$. Thus we   parameterize the
subgrid scale term $R(x,t)$, which is time-correlated, by colored
noise $\dot{B}_t^H$  as follows:
\begin{eqnarray}
\label{noise}
 R(x,t) = f(\bar{u}) + \sigma(x) \frac{d B^H_t}{dt},
\end{eqnarray}
  where
\begin{eqnarray} \label{force}
  f(\bar{u})=\EX R(x,t),
\end{eqnarray}
   is the mean component
  of the subgrid scale term $R(x,t)$. Moreover,
    the noise intensity $\sigma(x)$ is a  non-negative  deterministic function
to be determined from fluctuating SGS data $R$. The subgrid scale
term  $R(x, t)$ may be inferred from observational data
\cite{Peters1, Peters2}, or from fine mesh simulations as we do
here.  We represent the mean component $f(\bar{u})$ in terms of
the large eddy solution $\bar{u}$. The specific form for $f$
depends on the nature of the mean of $R$. Here we take
$f(\bar{u})= a_0+ a_1 u +a_2 u^2 +a_3 u^3$, where coefficients
$a_i$'s are determined via data fitting by minimizing $\int_0^T
\int_D[a_0+ a_1 u +a_2 u^2 +a_3 u^3 - \EX R(x,t)]^2dx dt$.
Moreover, we take $B_t^H$ as  a scalar fractional  Brownian
motion.
  % (i.e., spatially homogeneous).

  %With this  stochastically parameterized SGS term $R$, the large eddy
  %solution $\bar{u}$ from (\ref{LES}) is a stochastic process as well.

  %Note that  $\bar{u}$, as a solution of a \emph{nonlinear} random
  %partial differential equation, is generally non-Gaussian. The
  %product $\sigma(\bar{u}, x,t) \dot{B}^H_t$ is generally also
  %non-Gaussian. Overall, $\sigma(\bar{u}, x,t) \dot{B}^H_t$
  %may approximate quite complicated random data $R$.

Note that $\sigma$ is to be calculated or estimated from the
fluctuating SGS data $R$, either from observation or (in this
paper) from fine mesh simulations; see detailed discussions in
\cite{Menev, DuanBalu}. So this is an inverse problem. As in usual
inverse problems \cite{Tar}, the stochastic parameterizations for
the SGS term $R$ is not unique. This offers an opportunity for
trying various stochastic parameterization schemes, much as one
uses various smoother functions (e.g., polynomials or Fourier
series) to approximate less regular functions or data in
deterministic approximation theory.

To estimate the unknown parameter (function) $\sigma(x)$, we start
with  \eqref{noise}-\eqref{force} to get the following relation:
\begin{eqnarray}
 R(x,t) - \EX R(x,t) &=& \sigma(x) \frac{d B^H_t}{dt}.
\end{eqnarray}
 Taking  time integral over a computational interval $[0, T]$
  on  both  sides, we obtain
\begin{eqnarray*}
  \int_0^T[R(x,t)-\EX R(x,t)]dt  =  \int_0^T \sigma(x)
 d{B}^H_t = \sigma(x) B^H_T.
\end{eqnarray*}
Therefore, taking mean-square on both sides,
 \begin{eqnarray*}
  \EX( \int_0^T [R(x,t)-\EX R(x,t)]dt)^2 &=& \sigma^2(x)  T^{2H}.
\end{eqnarray*}
%Applying the Ito  isometry for stochastic integral with respect to
%fractional Brownian motion with $H > \frac12$, and using the fact
%that the integrand is deterministic, we obtain that
%\begin{eqnarray*}
 % \EX( \int_0^T [R(x,t)-\EX R(x,t)]dt)^2 &=&
% \sigma^2(x)\;H(2H-1)\int_0^T  \int_0^T |u-v|^{2H-2}dudv.
%\end{eqnarray*}
Thus an  estimator  for $\sigma(x)$ is
\begin{eqnarray}
 \sigma(x)= \frac1{T^H}\sqrt{\EX( \int_0^T [R(x,t)-\EX R(x,t)]dt)^2} \;\;,   \label{ax}
\end{eqnarray}
which can be computed numerically.

%With the calculated noise intensity term $\sigma(x)$ from
%(\ref{ax}), the stochastic LES model (\ref{new}) is established.

By the stochastic   parameterization  (\ref{noise}) on the SGS
term $R$, with  $f$ determined from (\ref{force}) and $\sigma $
from (\ref{ax}), the LES model (\ref{LES}) becomes a stochastic
partial differential equation (SPDE) for the large eddy solution
$U \approx \bar{u}  $:
\begin{eqnarray} \label{new}
  U_t = U_{xx} +  U - U^3
   + \int^t_0 {1 \over {1+|t-s|^{\beta}}}\; U(x,s)ds
   + f(U) + \sigma(x)   \frac{d B^H_t}{dt},
\end{eqnarray}
with boundary conditions $U(-1,t)=a, \; U(1,t)=b$  and filtered
initial condition
\begin{eqnarray} \label{newIC}
 U (x,0)= \bar{u}_0(x).
\end{eqnarray}

 {\bf Numerical Experiments:}

 We use a spectral method to solve nonlinear
 system (\ref{heat}) and (\ref{new})
  numerically.
For more details, please see \cite{Trefethen}. We take the
following initial and boundary conditions:
$$ u(x, 0)=u_0 = 0.53x
-0.47sin(1.5\pi x), \;\;
 u(-1,t)=-1, \;\;
 u(1,t)=1
$$

Fine mesh simulations of the original system with memory
(\ref{heat}) are conducted  to generate benchmark solutions or
solution realizations, with initial conditions slightly perturbed;
see Fig. 3. These fine mesh solutions $u$ are used to generate the
SGS term $R$ defined in (\ref{EF}) at each time and space step.
The filter size used in calculating $R$ is taken as $\delta=0.01$.
The mean $f$ is calculated from (\ref{force}) via cubic polynomial
data fitting (as discussed in the last section), and parameter
function $\sigma(x)$ is calculated as in (\ref{ax}). The
stochastic LES model (\ref{new}) is solved by the same numerical
code but on a coarser mesh.
%Various realizations of the
%solutions are obtained by running the same program again. The
%different simulation of fractional Brownian motion paths will give
%us different solution realizations for the LES model, since the
%solution also depends on the random paths.
Note that a four times coarser mesh simulation with no stochastic
parameterization for the original system (\ref{heat}) does not
generate satisfactory results; see Fig. 4. The stochastic LES
model    (\ref{new}) is then solved in the   mesh four times
coarser than the fine mesh used to solve the original
  equation (\ref{heat}). The stochastic parameterization leads to
  better resolution of the solution as shown in Fig. 5.
%  The root-mean-square error, $error(x, t):=\sqrt{\EX |u(x, t)-U(x,t)|^2}$,
%between the fine mesh solution $u$ (Fig.3)  and this stochastic
%LES solution $U$ (Fig. 5) is plotted in Fig. 6.
As in \cite{DuanBalu}, it can be shown that when two stochastic
parameterization terms are close in mean-square norm on finite
time intervals, the solutions are also close in the same norm.

\begin{figure}[htbp]
\begin{center}
\includegraphics[height=3in,width=4in]{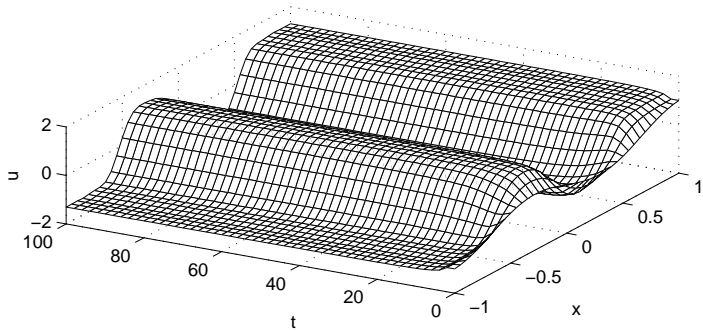}
\caption{\small \sl Solution to the original system on a fine
mesh, $u_t= u_{xx} + u - u^3 + \int^t_0 {1 \over
{1+|t-s|^{\beta}}} \;u(x,s)ds$, $\beta=2$,   mesh size $\Delta
x=0.001$. }
\end{center}
\end{figure}

\begin{figure}[htbp]
\begin{center}
\includegraphics[height=3in,width=4in]{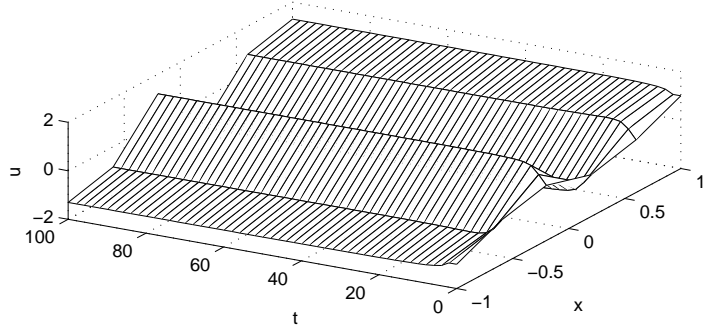}
\caption{\small \sl Solution to  the original system with NO
stochastic parametrization on the mesh four times coarser than the
mesh used in Fig. 3, $u_t = u_{xx} + u - u^3 + \int^t_0 {1 \over
{1+|t-s|^{\beta}}}\; u(x,s)ds$, $\beta=2$.}
\end{center}
\end{figure}

\begin{figure}[htbp]
\begin{center}
\includegraphics[height=3in,width=4in]{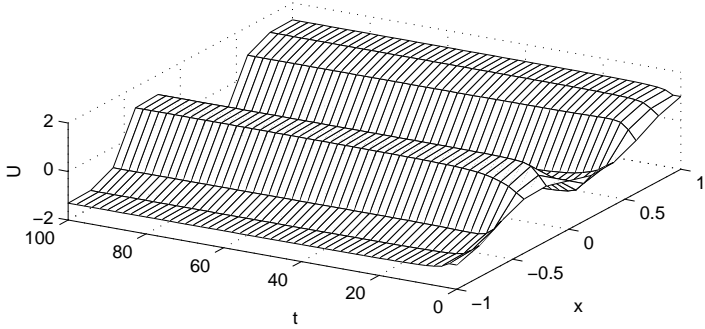}
\caption{\small \sl Solution to  LES model with stochastic
parametrization on the mesh four times coarser than the mesh used
in Fig. 3, $U_t = U_{xx} + U - U^3  + \int^t_0  {1 \over
{1+|t-s|^{\beta}}} \;U(x,s)ds + f(U)+ a(x) \dot{B}^{H}_t$,
$\beta=2$, $H=\frac34$.}
\end{center}
\end{figure}

%\begin{figure}[htbp]
%\begin{center}
%\includegraphics[height=3in,width=4in]{4X-coarse-error.eps}
%\caption{\small \sl Mean-square error of the LES model on the mesh
%four times coarser than the mesh used in Fig. 3, $U_t = U_{xx} + U
%- U^3  + \int^t_0  {1 \over {1+|t-s|^{\beta}}} \;U(x,s)ds +  f(U)
%+a(x) \dot{B}^{H}_t$, $\beta=2$, $H=\frac34$.}
%\end{center}
%\end{figure}

%\begin{figure}[htbp]
%\begin{center}
%\includegraphics[height=3in,width=4in]{old.eps}
%\caption{\small \sl Solution to the Ginzburg-Landau equation (\ref{heat}) on a fine mesh (dx=0.05)with $u(x, 0) = 0.53x
%-0.47sin(1.5\pi x), \;\;
% u(-1,t)=-1, \;\;
% u(1,t)=1.$}
%\end{center}
%\end{figure}

\bibliographystyle{amsplain}
%    Insert the bibliography data here.

\end{document}